\documentclass[10pt,a4paper]{article}

\usepackage[utf8]{inputenc}
\usepackage[T1]{fontenc}

\usepackage{bbm}
\usepackage{hyperref}
\usepackage{amsmath}
\usepackage{amsfonts}
\usepackage{amsthm}
\usepackage{amssymb}
\usepackage{graphicx}

\usepackage{bbm}
\newcommand{\R}{\mathbb R}
\newcommand{\ind}[1]{{\mathbbm{1}}_{#1}}
\newcommand{\bs}[1]{{\boldsymbol{#1}}}

\theoremstyle{plain}
\newtheorem{Theorem}{Theorem}
\newtheorem{Corollary}{Corollary}
\newtheorem{Proposition}{Proposition}

\theoremstyle{definition}
\newtheorem{Definition}{Definition}

\begin{document}

\title{The Rescaled P\'olya Urn and the Wright-Fisher process with mutation}

\author{Giacomo Aletti
\footnote{ADAMSS Center,
  Universit\`a degli Studi di Milano, Milan, Italy, giacomo.aletti@unimi.it}
$\quad$and$\quad$ 
Irene Crimaldi
\footnote{IMT School for Advanced Studies, Lucca, Italy, 
irene.crimaldi@imtlucca.it}
} 
\maketitle

\abstract{In \cite{ale-cri-RP,ale-cri-sar-sentiment} the authors introduce, study and apply a new
variant of the Eggenberger-P\'olya urn, called the ``Rescaled''
P\'olya urn, which, for a suitable choice of the model parameters, is
characterized by the following features: (i) a ``local''
reinforcement, i.e.~a reinforcement mechanism mainly based on the last
observations, (ii) a random persistent fluctuation of the predictive
mean, and (iii) a long-term almost sure convergence of the empirical
mean to a deterministic limit, together with a chi-squared goodness of
fit result for the limit probabilities. In this work, motivated by some empirical evidences in \cite{ale-cri-sar-sentiment}, we show that the
multidimensional Wright-Fisher diffusion with mutation can be obtained as a
suitable limit of the predictive means associated to a family of
rescaled P\'olya urns.}\\[10pt]
\noindent {\em Keywords:}
{P\'olya urn; predictive mean; urn model; Wright-Fisher diffusion.}

\section{Introduction}
\begin{figure}[tbp]
\begin{center}
\fbox{\includegraphics[width= 0.65\textwidth]{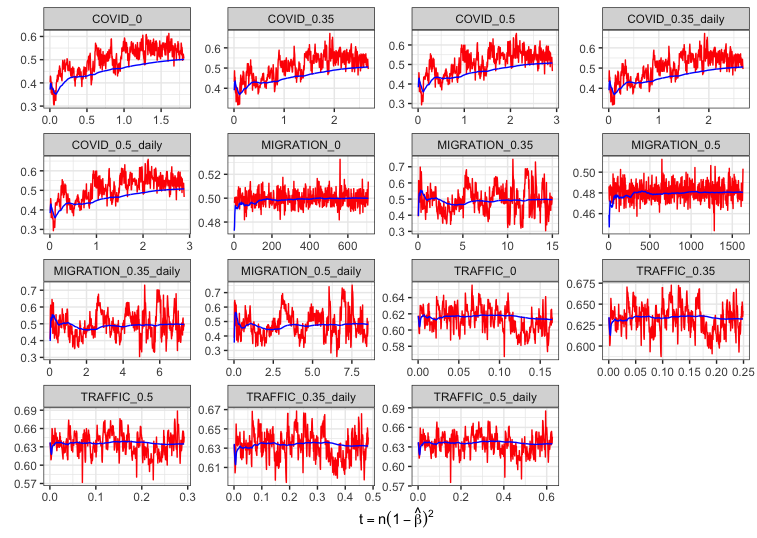}}
\end{center}
\caption{Twitter data: In \cite{ale-cri-sar-sentiment} the RP urn has been proven to be a good model for the evolution of the sentiment associated to Twitter posts. For these processes we have estimated   values of $\beta$ smaller than $1$, but very near to $1$. We here plot the processes $(\psi_{n\, 1})_n$ (red color) and $(\bar{\xi}_{n\,1} )_n$ (blue color), reconstructed from the data and rescaled in time  as $t= n (1-\beta)^2$. Details about the analyzed data sets, the reconstruction process and the estimated parameters can be found in \cite{ale-cri-sar-sentiment}.}
\label{figura}
\end{figure}

The standard Eggenberger-P\'olya urn \cite{EggPol23, mah} has
been widely studied and generalized. 
%
%
In its simplest form, this model with $k$-colors works as follows.
An urn contains $N_{0\, i}$ balls of color $i$, for $i=1,\dots, k$,
and, at each time-step, a ball is extracted from the urn and then
it is returned inside the urn together with $\alpha>0$ additional
balls of the same color. Therefore, if we denote by $N_{n\, i}$ the
number of balls of color $i$ in the urn at time-step $n$, we have
\begin{equation*}
N_{n\, i}=N_{n-1\,i}+\alpha\xi_{n\,i}\qquad\mbox{for } n\geq 1,
\end{equation*}
where $\xi_{n\,i}=1$ if the extracted ball at time-step $n$ is of color
$i$, and $\xi_{n\,i}=0$ otherwise. The parameter $\alpha$ regulates
the reinforcement mechanism: the greater $\alpha$, the greater the
dependence of $N_{n\,i}$ on $\sum_{h=1}^n\xi_{h\,i}$.\\
\indent In \cite{ale-cri-GRP, ale-cri-RP, ale-cri-sar-sentiment} the
Rescaled P\'olya (RP) urn has been introduced, studied, generalized
and applied.  This model is characterized by the introduction of a
parameter $\beta$ in the original model so that
\begin{equation*}
\begin{aligned}
N_{n\, i}& =b_{i}+B_{n\, i} &&\text{with }
\\
B_{n+1\, i}&=\beta B_{n\, i}+\alpha \xi_{n+1\, i}&& n\geq 0.
\end{aligned}
\end{equation*}
Therefore, the urn initially contains $b_{i}+B_{0\,i}>0$ balls of
color $i$ and the parameter $\beta\geq 0$, together with $\alpha>0$,
regulates the reinforcement mechanism. More precisely, the term $\beta
B_{n\,i}$ links $N_{n+1\,i}$ to the ``configuration'' at time-step $n$
through the ``scaling'' parameter $\beta$, and the term
$\alpha\xi_{n+1\,i}$ links $N_{n+1\,i}$ to the outcome of the
extraction at time-step $n+1$ through the parameter $\alpha$. The case
$\beta=1$ obviously corresponds to the standard Eggenberger-P\'olya
urn with an initial number $N_{0\,i}=b_{i}+B_{0\,i}$ of balls of color
$i$. When $\beta<1$, the RP urn model exhibits the following three
features:
\begin{itemize}
\item[(i)] a ``local'' reinforcement, i.e.~a reinforcement mechanism mainly based on the last
observations;
\item[(ii)] a random persistent fluctuation of the predictive
mean $\psi_{n\,i}=E[\xi_{n+1\,i}=1|\,\xi_{h\,j},\, 0\leq h\leq n,\,
  1\leq j\leq k]$; 
\item[(iii)] a long-term almost sure convergence of the empirical mean
  $\sum_{n=1}^N\xi_{n\, i}/N$ to the deterministic limit
  $p_{i}=b_{i}/\sum_{i=1}^nb_{i}$, and a chi-squared goodness of fit
  result for the long-term probability distribution
  $\{p_{1},\dots,p_{k}\}$.
\end{itemize}
Regarding point (iii), we specifically have that the chi-squared
statistics $$\chi^2 = N\sum_{i=1}^k \frac{(O_i/N-p_{i})^2}{p_{i}},$$
where $N$ is the size of the sample and $O_i=\sum_{n=1}^N\xi_{n\, i}$
the number of observations equal to $i$ in the sample, is
asymptotically distributed as $\chi^2(k-1)\lambda$, with $\lambda>1$.
This means that the presence of correlation among observations
mitigates the effect of the sample size $N$, that multiplies the
chi-squared distance between the observed frequencies and the expected
probabilities. This aspect is important for the statistical
applications in the context of a ``big sample'', when a small value of
the chi-squared distance might be significant, and hence a correction
related to the correlation between observations is desirable.  
%
%
In \cite{ale-cri-GRP, ale-cri-RP} it is described a possible
application in the context of clustered data, with independence
between clusters and correlation, due to a reinforcement mechanism,
inside each cluster.  \\
\indent In \cite{ale-cri-sar-sentiment} the RP urn has been applied as a good model for the evolution of the sentiment associated to Twitter posts. For these processes the estimated values of $\beta$ are strictly smaller than $1$, but very near to $1$. Note that the RP urn
dynamics with such a value for $\beta$ cannot be approximated by the
standard P\'olya urn ($\beta=1$), because one would loose the
fluctuations of the predictive means and the possibility of touching
the barriers $\{0,1\}$.  
%
%
In Figure \ref{figura}, we show the plots of the processes $(\psi_{n\, 1})_n$ and $(\bar{\xi}_{n\,1} )_n$, reconstructed from the data and rescaled in time  as $t= n (1-\beta)^2$. (Details about the analyzed data sets, the reconstruction process and the parameters estimation can be found in \cite{ale-cri-sar-sentiment}.) In this work, we show that the law of such processes can be approximated by the one of the Wright-Fisher diffusion with mutation. More precisely, we prove  
%
%
that the multidimensional Wright-Fisher
diffusion with mutation can be obtained as a suitable limit of the
predictive means associated to a family of RP urns with $\beta\in $[0,1)$, \beta\to 1$.  
\\
\indent The Wright–Fisher (WF) class of diffusion processes models
the evolution of the relative frequency of a genetic variant, or
allele, in a large randomly mating population with a finite number $k$
of genetic variants. When $k=2$, the WF diffusion obeys the
one-dimensional stochastic differential equation
\begin{equation}\label{WF-eq}
dX_t = F(X_t)dt + \sqrt{X_t(1-X_t)}dW_t,\qquad  X_0 = x_0, t \in [0,T].
\end{equation}
The drift coefficient, $F:[0, 1]\to R$, can include a variety of
evolutionary forces such as mutation and selection. For example,
$F(x)=p_1-(p_1+p_2)x=p_1(1-x)-p_2x$ describes a process with recurrent
mutation between the two alleles, governed by the mutation rates
$p_1>0$ and $p_2>0$. The drift vanishes when $x=p_1/(p_1+p_2)$ which
is an attracting point for the dynamics. Equation \eqref{WF-eq} can be
generalized to the case $k>2$.
%
%
The WF diffusion processes are widely employed in Bayesian Statistics,
as models for time-evolving priors \cite{fav-et-al, gri-spa, mena-rug,
  wal-et-al} and as a dicrete-time finite-population construction
method of the two-parameter Poisson-Dirichlet diffusion
\cite{cost-et-al}. They have been applied in genetics \cite{bol-et-al,
  gut-et-al, mala-et-al, sch-et-al, wil-et-al, zhao-et-al}, in
biophysics \cite{dan-et-al, dan-et-al-2}, in filtering theory
\cite{cha-gen, pap-rug} and in finance \cite{del-shi, gou-jas}.  \\
%
%
\indent The benefit coming from the proven limit result 
is twofold. First, the known properties of the WF process can
give a description of the RP urn when the parameter $\beta$ is
strictly smaller than one, but very near to one. Second, the given result might furnish the theoretical base for a new simulation method of the WF
process. Indeed, simulation from Equation \eqref{WF-eq} is highly
nontrivial because there is no known closed form expression for the
transition function of the diffusion, even in the simple case with
null drift \cite{jen-spa}.  \\
\indent The sequel of the paper is so structured. In Section
\ref{urn-model} we set up our notation and we formally define the RP
urn model. Section \ref{main} provides the main result of this work,
that is the convergence result of a suitable family of predictive
means associated to RP urns with $\beta\to 1$. In Section
\ref{properties} we list some properties of the considered stochastic
processes. In particular, we recall some properties of the WF
diffusion with mutation, connecting them to the parameters of the RP
urn model. Section \ref{case-two} focuses on the case $k=2$. Finally,
in Section \ref{dominant} we introduce the notion of dominant
component (color in the RP urn), related to the possibility of
reaching the barrier $1$. The paper closes with two technical
appendices.


\section{The Rescaled P\'olya urn}
\label{urn-model}

In all the sequel (unless otherwise specified) we suppose given two
parameters $\alpha>0$ and $\beta \geq 0$.  Given a vector $\bs{x}=
(x_1, \ldots, x_k)^\top\in \mathbb{R}^k$, we set \( |\bs{x}| =
\sum_{i=1}^k |x_i| \) and \( \|\bs{x}\|^2 = \bs{x}^\top \bs{x}=
\sum_{i=1}^k |x_i|^2 \).  Moreover we denote by $\bs{1}$ and $\bs{0}$
the vectors with all the components equal to $1$ and equal to $0$,
respectively. 
%
%
\\

To formally work with the RP urn model presented in the introduction,
we add here some notations.  In the whole sequel the expression
``number of balls'' is not to be understood literally, but all the
quantities are real numbers, not necessarily integers.  The urn
initially contains $b_{i}+B_{0\,i}>0$ distinct balls of color $i$,
with $i=1,\dots,k$.  We set $\bs{b}=(b_{1},\dots,b_{k})^{\top}$ and
$\bs{B_0}=(B_{0\,1},\dots,B_{0\,k})^{\top}$. In all the sequel (unless
otherwise specified) we assume $|\bs{b}|>0$ and we set $\bs{p} =
\frac{\bs{b}}{|\bs{b}|}$.  At each time-step $(n+1)\geq 1$, a ball is
drawn at random from the urn, obtaining the random vector
$\bs{\xi_{n+1}} = (\xi_{n+1\,1}, \ldots, \xi_{n+1\,k})^\top$ defined
as
\begin{equation*}
\xi_{n+1\,i} = 
\begin{cases}
1  &  \text{when the extracted ball at time $n+1$ is of color $i$}
\\
0  & \text{otherwise},
\end{cases}
\end{equation*}
and the number of balls in the urn is so updated:
\begin{equation}\label{eq:reinf1:K}
\bs{N_{n+1}}=\bs{b}+\bs{B_{n+1}}\qquad\text{with}
\qquad
\bs{B_{n+1}} = \beta \bs{B_n} + \alpha \bs{\xi_{n+1}}\,,
\end{equation}
which gives
\begin{equation}\label{eq:B}
\bs{B_{n}}=\beta^n \bs{B_0} + 
\alpha \beta^n \sum_{h=1}^n \beta^{-h} \bs{\xi_{h}}\,.
\end{equation}
Similarly, from the equality 
\begin{equation*}
|\bs{B_{n+1}}|= \beta |\bs{B_n}| + \alpha\,,  
\end{equation*}
we get, using $\sum_{h=0}^{n-1} x^h = (1-x^n)/(1-x)$, 
\begin{equation}\label{eq:|B|}
|\bs{B_n}|=\beta^n |\bs{B_0}| + \alpha \sum_{h=1}^n  \beta^{n-h}=
\beta^n\left(|\bs{B_0}|-\frac{\alpha}{1-\beta}\right)+\frac{\alpha}{1-\beta}\,.
\end{equation}
Setting $r^*_{n} = |\bs{N_n}|= |\bs{b}|+|\bs{B_n}|$, that
is the total number of balls in the urn at time-step $n$, we get the relations
\begin{equation}\label{dinamica-rnstar}
r_{n+1}^*=r_n^*+(\beta-1)|\bs{B_n}|+\alpha
\end{equation}
and
\begin{equation}\label{eq-rstar_n}
r_n^*=|\bs{b_0}|+\frac{\alpha }{1-\beta}+
\beta^n\left(|\bs{B_0}|- \frac{\alpha }{1-\beta}\right).
\end{equation}
Moreover, setting $\mathcal{F}_0$ equal to the trivial $\sigma$-field
and $\mathcal{F}_n=\sigma(\bs{\xi_1},\dots,\bs{\xi_n})$ for $n\geq 1$,
the conditional probabilities $\bs{\psi_{n}}= (\psi_{n\,1}, \ldots,
\psi_{n\,k})^\top$ of the extraction process, also called {\em predictive
means}, are
\begin{equation}\label{eq:extract1a:K-vettoriale}
\bs{\psi_{n}}=E[ \bs{\xi_{n+1}}| \mathcal{F}_n ]= 
\frac{\bs{N_n}}{|\bs{N_n}|} =
\frac{\bs{b}+\bs{B_n}}{r_n^*}\qquad n\geq 0
\end{equation}
and, from \eqref{eq:B} and \eqref{eq:|B|}, we have
\begin{equation}\label{eq-psi_n}
\bs{\psi_n} = 
\frac{ \bs{b} + \beta^n \bs{B_0} + \alpha \sum_{h=1}^n \beta^{n-h} \bs{\xi_{h}} }
     { |\bs{b}| + \frac{\alpha}{1-\beta} +
       \beta^n \left( |\bs{B_0}| - \frac{\alpha}{1-\beta} \right) }\,.
\end{equation}
The dependence of $\bs{\psi_n}$ on $\bs{\xi_h}$ depends on the factor
$f(h,n)=\alpha \beta^{n-h}$, with $1\leq h\leq n,\,n\geq 0$. In the
case of the standard Eggenberger-P\'olya urn, that corresponds to
$\beta=1$ for all $n$, each observation $\bs{\xi_h}$ has the same
``weight'' $f(h,n)=\alpha$. Instead, when $\beta<1$ the factor
$f(h,n)$ increases with $h$, then the main contribution is given by
the most recent extractions. We refer to this phenomenon as {\em
  ``local'' reinforcement}. The case $\beta=0$ is an extreme case, for
which $\bs{\psi_n}$ depends only on the last extraction $\bs{\xi_n}$.
\\ \indent By means of \eqref{eq:extract1a:K-vettoriale}, together
with \eqref{eq:reinf1:K} and \eqref{dinamica-rnstar}, we get 
\begin{equation}\label{dinamica-psin}
\bs{\psi_{n+1}}-\bs{\psi_{n}}=
-\frac{(1-\beta)}{r_{n+1}^*}|\bs{b}|\big(\bs{\psi_{n}}-\bs{p}\big)
+
\frac{\alpha}{r_{n+1}^*}\big(\bs{\xi_{n+1}}-\bs{\psi_{n}}\big).
\end{equation}
\indent Setting
%
%
$\Delta
\bs{M_{n+1}}= \bs{\xi_{n+1}}-\bs{\psi_{n}}$ 
%
%
and letting $\epsilon_n =|\bs{b}|(1-\beta)/r_{n+1}^*$ and $\delta_n
=\alpha/r_{n+1}^*$, from \eqref{dinamica-psin} we obtain
\begin{equation}\label{dinamica-psin-bis}
\bs{\psi_{n+1}}-\bs{\psi_{n}}=-\epsilon_{n}(\bs{\psi_{n}}-\bs{p})+
\delta_n\Delta\bs{M_{n+1}}\,.
\end{equation}


\section{Main result}\label{main}
Consider the RP urn with parameters $\alpha>0$, $\beta\in [0,1)$ and
  $\bs{B_0}$ such that $|\bs{B_0}| = r(\beta)=\alpha/(1-\beta)$ and
  set $b=|\bs{b}|>0$. Consequently, the total number of balls in the
  urn along the time-steps is constantly equal to
  $r^*(\beta)=b+r(\beta)$ and, if we denote by
  $\bs{\psi^{(\beta)}}=(\bs{\psi_n^{(\beta)}})_n$ the predictive means
  corresponding to the fixed value $\beta$, we have the dynamics
\begin{equation}\label{dinamica-psinEps}
\bs{\psi^{(\beta)}_{n}}-\bs{\psi^{(\beta)}_{n-1}}=
-\epsilon(\beta)\big(\bs{\psi^{(\beta)}_{n-1}}-\bs{p}\big)+
\delta(\beta)\Delta\bs{M_n^{(\beta)}}\,,
\end{equation}
where
\begin{equation}
\epsilon(\beta)=\frac{b(1-\beta)^2}{\alpha+b(1-\beta)},\qquad
\delta(\beta)=\frac{\alpha(1-\beta)}{\alpha+b(1-\beta)}
\end{equation}
and $\bs{\Delta
  M_n^{(\beta)}}=\bs{\xi^{(\beta)}_{n}}-\bs{\psi^{(\beta)}_{n-1}}$.
(Note that we have $\epsilon(\beta)\sim c \delta(\beta)^2$ for
$\beta\to 1$, with $c=b/\alpha>0$.)  Finally, define ${\bs
  X^{(\beta)}}=(\bs{X_t^{(\beta)}})_{t\geq 0}$ where
\begin{equation}
\bs{X^{(\beta)}_t} = \bs{\psi^{(\beta)}_{\lfloor t/(1-\beta)^2 \rfloor}} 
\quad \iff \quad 
\bs{X^{(\beta)}_t} = \bs{\psi^{(\beta)}_{n-1}}, \
t \in [\, (n-1) (1-\beta)^2, n(1-\beta)^2\,).
\end{equation}

The following result holds true:
\begin{Theorem}\label{teo:converg}
  Suppose that $\bs{X^{(\beta)}_0}$ weakly converges towards some
  process $\bs{X_0}$ when $\beta\to 1$. Then, for $\beta\to 1$, the
  family of stochastic processes $\{\bs{X^{(\beta)}}, \, \beta\in
  [0,1)\}$ weakly converges towards the $k$-alleles Wright-Fisher
    diffusion $\bs{X}=(\bs{X_t})_{t\geq 0}$, with 
    \emph{type-independent mutation kernel} given by $\bs{p}$
    %
    %
    and dynamics 
 \begin{equation}\label{eq:WLimit}
   d \bs{X_t} = - b\frac{\bs{X_t}  - \bs{p}}{\alpha} dt +
   \Sigma( \bs{X_t} ) d \bs{W_t} ,  
\end{equation}
with $\Sigma( \bs{X_t} ) \Sigma( \bs{X_t} )^{\top} = \Big(
\mathrm{diag} (\bs{X_t} ) - \bs{X_t} \bs{X_t}^{\top}\Big)$ and
$\bs{1}^\top \Sigma( \bs{X_t} ) = \bs{0}^\top$, that is
\begin{align}\label{eq:matrix}
& \Sigma( \bs{X_t} )_{ij} = 
\begin{cases}
0
  & \text{if $X_{t,i} X_{t,j} = 0$ or $i<j$}
\\
\sqrt{ X_{t,i} \frac{\sum_{l = i+1}^k X_{t,l} }{ \sum_{l = i}^k X_{t,l} } }
& \text{if $i=j$ and $X_{t,i} X_{t,j} \neq 0$}
\\
- X_{t,i} \sqrt{ \frac{X_{t,j} }{ \sum_{l = j}^k X_{t,l} \sum_{l=j+1}^k X_{t,l} } }
& \text{if $i>j$ and $X_{t,i} X_{t,j} \neq 0$.}
\end{cases}
\end{align}
%
%
\end{Theorem}

\begin{proof}
  Fix a sequence $(\beta_n)$, with $\beta_n\in[0,1)$ and $\beta_n \to
    1$.  The sequence of processes $\{ \bs{X^{(\beta_n )}},\, n \in
    \mathbb{N}\}$ is bounded, and hence we have to prove the tighthness
    of the sequence in the space $D^k[0,\infty)$ of right-continuous
      functions with the ususal Skorohod topology, and the
      characterization of the law of the unique limit process.

For any $f\in C^2_b$, define
\begin{equation}\label{eq:generators}
\begin{aligned}
\gamma_n^{(\beta,f)}(\bs{x}) & = 
\widehat{A}^{(\beta)}f((n-1) (1-\beta)^2) (\bs{x})
\\
& = E\Big[ 
  \frac{f(\bs{X^{(\beta)}_{n(1-\beta)^2}} )-
    f(\bs{X^{(\beta)}_{(n-1) (1-\beta)^2}}) }{(1-\beta)^2} 
\Big| \bs{X^{(\beta)}_{(n-1) (1-\beta)^2}}  = \bs{x}\Big] \\
& =
E\Big[ 
\frac{f(\bs{\psi^{(\beta)}_{n}})-f(\bs{\psi^{(\beta)}_{n-1}}) }{(1-\beta)^2} 
\Big| \bs{\psi^{(\beta)}_{n-1}} = \bs{x}\Big] \\
& \mathop{=}^{\text{by \eqref{dinamica-psinEps}}} 
\frac{1}{(1-\beta)^2} \bigg( E \Big[ f(\bs{x})  + 
\sum_i \frac{\partial f }{\partial x_i} (\bs{x}) ( -\epsilon(\beta)(x_i-{p_{i}} )
+ \delta(\beta) \Delta {M_{n,i}}^{(\beta)} ) 
\\
& \qquad \qquad 
+ \tfrac{1}{2}\delta(\beta)^2
\sum_{ij} \frac{\partial^2 f }{\partial x_i\partial x_j} (\bs{x})
\Delta M_{n,i}^{(\beta)} 
\Delta M_{n,j}^{(\beta)} 
+ O((1-\beta)^{3}) \Big|\mathcal{F}_{n-1} \Big] - f(\bs{x}) \bigg)
\\ &
= 
-\tfrac{b}{\alpha+b(1-\beta)}
\sum_i \frac{\partial f }{\partial x_i} (\bs{x}) ( x_i-{p_{i}} )
+ \frac{1}{2} \tfrac{\alpha^2}{(\alpha+b(1-\beta))^2}
\sum_{ij}
\frac{\partial^2 f }{\partial x_i\partial x_j} (\bs{x})( x_i\ind{i=j} - x_ix_j)\\
&+ O(1-\beta)\,.
\end{aligned}
\end{equation}
We note that, for any $f\in C^2_b$, the partial derivatives in
\eqref{eq:generators} are uniformly dounded, as $\bs{x}$ belongs to
the compact simplex $S = \{x_i\geq 0,\sum_i x_i = 1\}$. The family
$\{\gamma_n^{(\beta,f)}(\bs{x}) , n \in \mathbb{N}, \beta<1, \bs{x}
\in S \}$ is then uniformly integrable.  Thus, as a consequence of
\cite[Theorem 4]{Kus} (or \cite[ch. 7.4.3, Theorem~4.3,
  p.~236]{KusYin}), we have that the sequence of processes
$\{\bs{X^{(\beta_n )}} , n \in \mathbb{N}\}$ is tight in the space of
right-continuous functions with the ususal Skorohod topology.  Since,
for any $n$ and $t$, $\bs{X^{(\beta_n )}_t} \in S$, then $\bs{1}^\top
\Sigma( \bs{X_t} ) = \bs{0}^\top $.  Moreover, the generator of the
limit process is determined by the limit
\[
\begin{aligned}
  Af(t)(\bs{x}) & =
  \lim_{n\to\infty} \gamma_{\lfloor t/(1-\beta)^2 \rfloor}^{(\beta_n,f)}(\bs{x})  
\\
& = 
-\tfrac{b}{\alpha}
\sum_i \frac{\partial f }{\partial x_i} (\bs{x}) ( x_i-{p_{i}} )
+
\frac{1}{2}
\sum_{ij} \frac{\partial^2 f }{\partial x_i\partial x_j}
(\bs{x}) ( x_i\ind{i=j} - x_ix_j) .
\end{aligned}
\]
Hence, the weak limit of the sequence of the bounded processes 
$\bs{X^{(\beta_n )}}$ is the diffusion process
\[
d \bs{X_t} =
- b\frac{\bs{X_t}  - \bs{p}}{\alpha} dt + \Sigma( \bs{X_t} ) d \bs{W_t},  
\qquad  \Sigma( \bs{X_t} ) \Sigma( \bs{X_t} )^{\top} = 
\Big( \mathrm{diag} (\bs{X_t} ) - \bs{X_t} \bs{X_t}^{\top}\Big).
\]
The expression \eqref{eq:matrix} follows from \cite[Corollary~3]{Tanabe}. 
\end{proof}


\section{Some properties}\label{properties}
 We list here some properties.
\subsection{Projections}

Let $\bs{J}= \{J_1,\ldots,J_{k_J}\}$, be a partition of
$\{1,\ldots,k\}$, in that $J_l \neq \varnothing$, $J_{i_1} \cap
J_{i_2} = \varnothing$, and $ \cup_{l=1}^{k_J} = \{1,\ldots,k\}$.
Here $k_j$ denotes the cardinality of $\bs{J}$ .  Define the $
k_{J}$-dimesional objects $(\bs{\psi^{(\beta,J)}_{n}})_n $,
$(\bs{\xi^{(\beta,J)}_{n}})_n $ and $\bs{p^{(J)}}$ as
\[
\left.
\begin{aligned}
{\psi^{(\beta,\bs{J})}_{n,i}} & =
\sum_{l\in J_i} {\psi^{(\beta)}_{n,l}} 
\\
{\xi^{(\beta,\bs{J})}_{n,i}}
& =
\sum_{l\in J_l} {\xi^{(\varepsilon)}_{n,l}} 
\\
{p_{i}^{(\bs{J})}}
& =
\sum_{l\in J_i} {p_{l}} 
\end{aligned}
\right\}
\qquad \text{for }
i = 1,\ldots, k_{J},
\]
and $\bs{X^{(\beta,\bs{J})}_t} = \bs{\psi^{(\beta,\bs{J})}_{\lfloor
    t/(1-\beta)^2 \rfloor}} $.  With these definitions, from
\eqref{dinamica-psinEps}, we immediately get that
$(\bs{\psi^{(\beta,\bs{J})}_{n}})_n$ is a $k_J$-dimensional RP urn
following the dynamics
\begin{equation}\label{dinamica-projection}
\bs{\psi^{(\beta,\bs{J})}_{n}}-\bs{\psi^{(\beta,\bs{J})}_{n-1}}=
-\epsilon(\beta)\big(\bs{\psi^{(\beta,\bs{J})}_{n-1}}-\bs{p^{(J)}}\big)
+
\delta(\beta)\big(\bs{\xi^{(\beta,\bs{J})}_{n}}-\bs{\psi^{(\beta,\bs{J})}_{n-1}}\big)
\end{equation}
and that Theorem~\ref{teo:converg} holds for
$\bs{X^{(\beta,\bs{J})}_t}$.  Consequently, the convergence to the
Wright-Fisher diffusion still holds if we group together some
components of the process.  This property is summarized in the
following theorem:
\begin{Theorem}\label{teo:grouping}
Under the hypothesis of Theorem~\ref{teo:converg}, the process
$\bs{X^{(\beta,\bs{J})}_t}$ weakly converges to $\bs{X^{(\bs{J})}_t}$ which satisfy
the SDE
 \begin{equation}\label{eq:WLimitGroup}
   d \bs{X^{(\bs{J})}_t} = - b\frac{\bs{X^{(\bs{J})}_t}  - \bs{p^{(J)}}}{\alpha} dt +
   \Sigma( \bs{X^{(\bs{J})}_t} ) d \bs{W^*_t} ,
\end{equation}
where $\bs{W^*_t} $ is a $k_{J}$-dimensional standard Brownian motion, 
and $X^{(\bs{J})}_{t,i} = \sum_{l\in J_i} X_{t,l}$.
\end{Theorem}

\subsection{Limiting ergodic distribution}
Since the simplex has dimension $k-1$ with respect to the Lebesgue
measure, it is convenient to change the notations. Let $T^{k-1}$ be
the $k-1$-dimensional simplex defined by
$$
T^{k-1} :=
\{\bs{y} \in \R^{k-1}:
y_1 \ge 0, \ldots, y_{k-1} \ge 0, 1-y_1-y_2-\cdots-y_{k-1} \ge 0 \},  
$$ where, with the old definition, we have $x_i = y_i, i < k$ and $x_k
:= 1-y_1-y_2-\cdots-y_{k-1} $.  Obviously, there is a one-to-one
natural correspondence between $ T^{k-1}$ and the simplex $\{\bs{x}
\in \R^{k}: x_1 \ge 0, \ldots, x_{k} \ge 0, \sum_i x_i = 1 \}$ defined
by
\[
\bs{y} = (y_1 , \ldots, y_{k-1} ) \quad \longleftrightarrow \quad
(y_1 , \ldots, y_{k-1} , 1-y_1-y_2-\cdots-y_{k-1} ) = 
(x_1 , \ldots, x_{k-1}, x_{k}) = \bs{x}. 
\]
The Markov diffusion process $\bs{X_t}$ in \eqref{eq:WLimit} may be
ridefined as $\bs{Y_t}= (X_{t,1},\ldots,X_{t,k-1})$ on $\bs{y} \in
T^{k-1}$ with the corresponding generator
 \begin{equation}\label{eq:LimitGenerator}
Lf(\bs{y}) = 
-\tfrac{b}{\alpha} \sum_{i=1}^{k-1} \frac{\partial f }{\partial y_i}
(\bs{y}) ( y_i-{p_{i}} )
+ \frac{1}{2} \sum_{i,j=1}^{k-1} \frac{\partial^2 f }{\partial y_i\partial y_j}
(\bs{y}) ( y_i\ind{i=j} - y_iy_j) .
\end{equation}
The Kolmogorov forward equation for the density $p(\bs{y},t)$ of the
limiting process $\bs{Y}_t$ is
\begin{multline}\label{eq:KElimit}
  \frac{\partial}{\partial t} p(\bs{y},t) =
\frac{1}{2} \bigg(
\tfrac{b}{\alpha} \sum_{i=1}^{k-1} \frac{\partial  }{\partial y_i}
\Big( p(\bs{y},t) ( y_i-{p_{i}} ) \Big)
\\
+ \sum_{i=1}^{k-1} \frac{\partial^2 }{\partial y_i^2}
\Big(y_i (1 - y_i) p(\bs{y},t) \Big) 
- 2\sum_{1 \leq i < j \leq k-1} \frac{\partial^2 }{\partial y_i\partial y_j}
\Big(y_iy_j p(\bs{y},t) \Big) \bigg).
\end{multline}
Therefore, it is not hard to show that the limit invariant ergodic
distribution is
\begin{equation}\label{eq:limit-distr}
p(\bs{y}) =
\frac{1}{B(2\tfrac{b}{\alpha}\bs{p})}
(1-y_1-\cdots-y_{k-1})^{\frac{2 b{ (1- p_{1} - \cdots - p_{k-1} ) }}{\alpha} -1 }
\prod_{i=1}^{k-1} y_i^{\frac{2 b{p_{i}}}{\alpha} -1 } ,
\end{equation}
because it satisfy \eqref{eq:KElimit} (see also \cite{wright84}). The
above distribution is the Dirichel distribution $\hbox{Dir}\big( 2
\frac{b}{\alpha} \bs{p} \big)$ as a function of
$\bs{x}=(\bs{y},1-y_1-\cdots-y_{k-1})$.

%
%


\subsection{Transition density of the limit process}
The transition density $p(\bs{y_0},\bs{y};t)$ is defined by
\[
P(\bs{Y_t} \in S | \bs{Y_0}=\bs{y_0}) =
\int_{S\cap T^{k-1}} p(\bs{y_0},\bs{y};t) d\bs{y}
\]
and it can be represented in terms of series of orthogonal
polynomials given in Appendix~\ref{sec:polynomials}.  We first note that the
limiting invariant ergoding distribution $p(\bs{y})$ in
\eqref{eq:limit-distr} and the generator of the process $\bs{Y_t}$ in
\eqref{eq:LimitGenerator} may be rewritten on $T^{k-1}$ in terms of
$\gamma_i = 2\tfrac{b}{\alpha} p_i -1$, obtaining
 \begin{equation*}
\begin{aligned}
\pi_{\boldsymbol{\gamma}}(\bs{{y}}) & = 
\frac{1}{B(\boldsymbol{\gamma}+1)} 
\Big( \prod_{i=1}^{k-1} y_i^{\gamma_i}  \Big) (1- y_1-y_2-\cdots-y_{k-1} )^{\gamma_{k}} 
\\
Lf(\bs{y})
& = 
-\tfrac{b}{\alpha} \sum_{i=1}^{k-1} \frac{\partial f }{\partial y_i}
(\bs{y}) ( y_i-{p_{i}} )
+ \frac{1}{2} \sum_{i,j=1}^{k-1}
\frac{\partial^2 f }{\partial y_i\partial y_j} (\bs{y}) ( y_i\ind{i=j} - y_iy_j) 
\\
& = \frac{1}{2} \bigg(
\sum_{i=1}^{k-1} \Big( 2\tfrac{b}{\alpha} p_i -
\Big[\big(\sum_{i=1}^k 2\tfrac{b}{\alpha} {p_{i}} - 1\big) + k \Big] y_i \Big) 
\frac{\partial f }{\partial y_i} (\bs{y} )
\\
& \qquad\qquad
+ \sum_{i=1}^{k-1} y_i (1 - y_i) \frac{\partial^2 f }{\partial y_i^2} (\bs{y}) 
- 2\sum_{1 \leq i < j \leq k-1} y_iy_j
\frac{\partial^2 f }{\partial y_i\partial y_j} (\bs{y}) \bigg) 
\\
& = \frac{1}{2} \bigg(
\sum_{i=1}^{k-1} \Big( \gamma_i +1 - \Big[ \big(\sum_{i=1}^k \gamma_i \big)
  + k \Big] y_i \Big) 
\frac{\partial f }{\partial y_i} (\bs{y} )
\\
& \qquad\qquad
+ \sum_{i=1}^{k-1} y_i (1 - y_i) \frac{\partial^2 f }{\partial y_i^2} (\bs{y}) 
- 2\sum_{1 \leq i < j \leq k-1} y_iy_j
\frac{\partial^2 f }{\partial y_i\partial y_j} (\bs{y}) \bigg) .
\end{aligned}
 \end{equation*}
 These two expressions coincide with those given in
 %
 %
 Appendix~\ref{sec:polynomials}. Let
$\mathcal{V}_{n,\boldsymbol{\gamma}}$ be the space of orthogonal
polynomials of degree $n$ as defined there and let
$f^{\boldsymbol{\gamma}}_{{\mathbf{n}}} $ one of the three orthogonal
bases given there.  Then \eqref{eq:diff-eqn2} implies 
\[
L f^{\boldsymbol{\gamma}}_{{\mathbf{n}}} =
\frac{1}{2} L_{\boldsymbol{\gamma}} f^{\boldsymbol{\gamma}}_{{\mathbf{n}}} 
= \frac{1}{2}  (-\lambda_n) f^{\boldsymbol{\gamma}}_{{\mathbf{n}}}  =
-\nu_n f^{\boldsymbol{\gamma}}_{{\mathbf{n}}} ,
\]
where $ \nu_n = \frac{n(n+ k + \sum_{i=1}^k \gamma_i )}{2} =
\frac{n(n+2\tfrac{b}{\alpha}-1)}{2}$.  Note that each
$\psi^{\boldsymbol{\gamma}}_{{\mathbf{n}}} (t,\bs{y}) = e^{-\nu_n t}
f^{\boldsymbol{\gamma}}_{{\mathbf{n}}} (\bs{y}) $ satisfies the
Kolmogorov backward equation associated to the process $\bs{Y_t}$,
since
\[
\frac{\partial}{\partial t} \psi^{\boldsymbol{\gamma}}_{{\mathbf{n}}} (t,\bs{y}) = 
- \nu_n \psi^{\boldsymbol{\gamma}}_{{\mathbf{n}}} (t,\bs{y}) =
e^{-\nu_n t}  L f^{\boldsymbol{\gamma}}_{{\mathbf{n}}} (\bs{y})  =
L \psi^{\boldsymbol{\gamma}}_{{\mathbf{n}}} (t,\bs{y}) .
\]
Now, for any $\bs{y}\in T^{k-1}$, let $S_{\bs{y}} =
\{\bs{\bar{y}}\colon 0\leq \bar{y}_i \leq y_i, i= 1,\ldots,k-1\}$. The
function
\[
U_{\bs{y}}(\bs{y_0},t) = P(\bs{Y_t} \in {S_{\bs{y}}} | \bs{Y_0}=\bs{y_0}) =
\int_{T^{k-1}} \ind{S_{\bs{y}}}(\bs{\bar{y}}) 
p(\bs{y_0},\bs{\bar{y}};t) d\bs{\bar{y}}
\]
satisfies the Kolmogorov backward equation. As in
\cite[Section~15.13]{KarTay81}, we formally write
$U_{\bs{y}}(\bs{y_0},t) $ in terms of
$\psi^{\boldsymbol{\gamma}}_{{\mathbf{n}}} (t,\bs{y_0}) $, since the
Kolmogorov backward equation is additive, obtaining
\[
U_{\bs{y}}(\bs{y_0},t) = \sum_n \sum_{{\mathbf{n}}\colon n_1+\cdots+n_{k-1}=n} 
c_{\mathbf{n}}(\bs{y}) \psi^{\boldsymbol{\gamma}}_{{\mathbf{n}}} (t,\bs{y_0}) .
\]
Note that the boundary conditions imply that
\[
\begin{aligned}
 \ind{S_{\bs{y}}}(\bs{{y_0}}) 
& = 
\lim_{t\to 0^+} U_{\bs{y}}(\bs{y_0},t)  = 
\lim_{t\to 0^+} \sum_n \sum_{{\mathbf{n}}\colon n_1+\cdots+n_{k-1}=n} 
c_{\mathbf{n}}(\bs{y}) e^{-\nu_n t}  f^{\boldsymbol{\gamma}}_{{\mathbf{n}}} (\bs{y_0}) 
\\
& =
\sum_n \sum_{{\mathbf{n}}\colon n_1+\cdots+n_{k-1}=n} 
c_{\mathbf{n}}(\bs{y}) f^{\boldsymbol{\gamma}}_{{\mathbf{n}}} (\bs{y_0}) .
\end{aligned}
 \]
The orthogonality and the completeness of the polynomial system implies that
\[
c_{\mathbf{n}}(\bs{y}) = \frac{
\langle 
\ind{S_{\bs{y}}} ,
f^{\boldsymbol{\gamma}}_{{\mathbf{n}}} 
\rangle_{\boldsymbol{\gamma}}
}{
\langle 
f^{\boldsymbol{\gamma}}_{{\mathbf{n}}} , f^{\boldsymbol{\gamma}}_{{\mathbf{n}}} 
\rangle_{\boldsymbol{\gamma}}
}
=
\frac{
\int_{S_{\bs{y}}} 
f^{\boldsymbol{\gamma}}_{{\mathbf{n}}} (\bs{\bar{y}})
\pi_{\boldsymbol{\gamma}}(\bs{\bar{y}})
d\bs{\bar{y}}
}{
\langle 
f^{\boldsymbol{\gamma}}_{{\mathbf{n}}} , f^{\boldsymbol{\gamma}}_{{\mathbf{n}}} 
\rangle_{\boldsymbol{\gamma}}
}.
\]
The transition density $p( \bs{y_0}, \bs{y} ; t )$ may be then
computed differentiating $U_{\bs{y}}(\bs{y_0},t)$, obtaining
(cfr.\ \cite[Eq.~(15.13.11)]{KarTay81})
\[
\begin{aligned}
p(\bs{y_0},\bs{y};t) & 
= \frac{\partial}{\partial\bs{y}} \sum_n \sum_{{\mathbf{n}}\colon n_1+\cdots+n_{k-1}=n} 
c_{\mathbf{n}}(\bs{y}) 
\psi^{\boldsymbol{\gamma}}_{{\mathbf{n}}} (t,\bs{y_0}) 
\\
& =  
\pi_{\boldsymbol{\gamma}}(\bs{y})
\sum_n 
e^{-\nu_n t} 
\sum_{{\mathbf{n}}\colon n_1+\cdots+n_{k-1}=n} 
\frac{
  f^{\boldsymbol{\gamma}}_{{\mathbf{n}}} (\bs{y})
  f^{\boldsymbol{\gamma}}_{{\mathbf{n}}} (\bs{y_0})  
}{
\langle 
f^{\boldsymbol{\gamma}}_{{\mathbf{n}}} , f^{\boldsymbol{\gamma}}_{{\mathbf{n}}} 
\rangle_{\boldsymbol{\gamma}}
}.
\end{aligned}
\]

\section{Two-dimensional urn}\label{case-two}
%
In the next proposition we point out the behavior obtained when we
look at the aggregated evolution of two groups $J_1,\,J_2$ of urn 
colors.
\begin{Proposition}\label{prop:2dimURN}
Let $\bs{J} = \{J_1,J_2\} $ with $J_1\neq \varnothing$, $J_2\neq
\varnothing$ and $J_1\cup J_2 = \{1,\ldots,k\}$. Under the hypothesis
of Theorem~\ref{teo:converg}, each component of the sequence of
processes $\bs{X^{(\beta,J)}_t}$ converges, for $\beta\to 1$, to the
one dimensional diffusion process with values in
$[0,1]$ that satisfies the SDE
\[
d{X}^{(J_1,J_2)}_{t,i} =
-b\frac{{X}^{(J_1,J_2)}_{t,i} -
  {p_{i}}}{\alpha} dt + (-1)^{i+1}
  \sqrt{\max(0,{X}^{(J_1,J_2)}_{t,i} (1-{X}^{(J_1,J_2)}_{t,i}))} d {W}_t.
\]
In addition, ${X}^{(J_1,J_2)}_{t,1} = \sum_{l\in J_1} X_{t,l}$ and 
${X}^{(J_1,J_2)}_{t,2} = \sum_{l\in J_2} X_{t,l}$.
\end{Proposition}
\begin{proof}
It is sufficient to apply Theorem~\ref{teo:grouping} and note 
that the process $\bs{X^{(\bs{J})}_t} $ satisfies the SDE
\eqref{eq:WLimitGroup}, that now reads
\[
d 
\Big(
\begin{smallmatrix}
X^{(\bs{J})}_{t,1}
\\
X^{(\bs{J})}_{t,2}
\end{smallmatrix}
\Big)
= 
-\frac{b}{\alpha} 
\Big(
\begin{smallmatrix}
X^{(\bs{J})}_{t,1} - {p_{1}^{(\bs{J})}}
\\
X^{(\bs{J})}_{t,2} - {p_{2}^{(\bs{J})}}
\end{smallmatrix}
\Big)
dt 
+ 
\Big(
\begin{smallmatrix}
\sqrt{X^{(\bs{J})}_{t,1} X^{(\bs{J})}_{t,2}}
&  0
\\
- \sqrt{X^{(\bs{J})}_{t,1} X^{(\bs{J})}_{t,2}}
&
0
\end{smallmatrix}
\Big)
d
\Big(
\begin{smallmatrix}
W_{t,1}
\\
W_{t,2}
\end{smallmatrix}
\Big). \qedhere
\]
\end{proof}

Now, if we further specialize the grouping choice to 
$\bs{J}=(\{i\}, \{1,\ldots,i-1,i+1,\ldots,k\})$, we get
\[
\begin{aligned}
{\psi^{(\beta,\bs{J})}_{n,1}} 
&=
{\psi^{(\beta)}_{n,i}} , 
&&& 
{\psi^{(\beta,\bs{J})}_{n,2}} 
&=
\sum_{l\neq i} {\psi^{(\beta)}_{n,l}} 
\\
{\xi^{(\beta,\bs{J})}_{n,1}} 
&=
{\xi^{(\beta)}_{n,i}} , 
&&& 
{\xi^{(\beta,\bs{J})}_{n,2}} 
&=
\sum_{l\neq i} {\xi^{(\beta)}_{n,l}} = 1 - 
{\xi^{(\beta)}_{n,i}}
\\
{p_{1}^{(\bs{J})}}
& =
{p_{i}}, 
&&& 
{p_{2}^{(\bs{J})}}
& =
\sum_{l\neq i} p_l = 
1 -{p_{i}}.
\end{aligned}
\]
We note that in this case ${Y}_{t,i} = X^{(\bs{J})}_{t,1} $ is the
first component of $\bs{X^{(\bs{J})}_t} $ and we get the following
corollary:

\begin{Corollary}\label{cor:margin}
Under the conditions of Theorem~\ref{teo:converg} the $i$-th component
of the sequence of processes $\bs{X^{(\beta)}}$ converges, for
$\beta\to 1$, to the one dimensional diffusion $(X_{t,i})_{t\geq 0}$
with values in $[0,1]$ satisfying the SDE
\[
d{X}_{t,i} =
-b\frac{{X}_{t,i} -
  {p_{i}}}{\alpha} dt + 
  \sqrt{\max(0,{X}_{t,i}(1-{X}_{t,i}))} d {W}_t.
\]
\end{Corollary}


\subsection{Excursions from $z_0$}\label{sec:group}
Let $\bs{J} = \{J_1,J_2\} $ as in Proposition~\ref{prop:2dimURN}, that
implies that $Z_t = \sum_{l\in J} X_{t,l}$
satisfies the following equation
\[
\begin{aligned}
d{Z}_{t} 
& =
-b\frac{{Z}_{t} -
  \sum_{l\in J}{p_{l}}}{\alpha} dt + 
  \sqrt{\max(0,{Z}_{t}(1-{Z}_{t}))} d {W}_t,
  .
  \end{aligned}
\]
that is \eqref{eq:procBarrier} with $a_0 = \frac{b}{\alpha} \sum_{l\in
  J}{p_{l}} $ and $ a_1 = \frac{b}{\alpha} - a_0 $.  We focus here on
some properties of $Z$. As in \cite[Section~15.3]{KarTay81}, let $a$
and $b$ be fixed, subject to $0< a < b < 1$, and let $\tau_A$ be the
hitting time of the set $A$ (for $z\in (0,1)$, we set $\tau_{z} =
\tau_{\{z\}}$ for semplicity) and $\tau^* = \tau_{\{a,b\}} =
\min(\tau_{a},\tau_{b})$ be the first time the process reaches either
$a$ or $b$.  We highlight some classical problems that are linked to
$\tau_b , \tau_a$ and $\tau^*$.

\subsubsection*{Problem 1.} Find $ u (z_0) = P (\tau_b < \tau_a |Z_0 = z_0)$, 
the probability that the process $Z_t$ reaches $b$ before $a$ starting
from $z_0$.  As a consequence of \cite[Eq.~(15.3.10)]{KarTay81} and of
the scale function $S:(0,1)\to\R$ of $Z_t$ given in \eqref{scale}, we
get
\[
u(z_0) = 
\frac{S(z_0)-S(a)}{S(b)-S(a)} 
=
\frac{
  \int_{a}^{z_0} t^{-2 \frac{b}{\alpha} \sum_{l\in J}{p_{l}} }
  (1-t)^{-2 \frac{b}{\alpha} (1-\sum_{l\in J}{p_{l}} )} dt 
}{
  \int_{a}^{b} t^{-2 \frac{b}{\alpha} \sum_{l\in J}{p_{l}} }
  (1-t)^{-2 \frac{b}{\alpha} (1-\sum_{l\in J}{p_{l}} )} dt 
} .
\]

\subsubsection*{Problem 2.} Find $ w(z_0) =
E ( \int_0^{\tau^*} g(Z_t) dt | Z_0 = z_0) $, with $g$ bounded 
and continuous function.  This quantity is the expected cost up to the
time when either $a$ or $b$ was first reached, under the cost rate
$g$, starting from $z_0 \in (a,b)$.  When $g\equiv 1$, then this
problem gives the mean time to reach either $a$ or $b$, starting from
$z_0 \in (a,b)$.  The solution $w:(a,b)\to \R$ may be computed in
terms of $u:(a,b)\to [0,1]$ above, of the scale function
$S:(0,1)\to\R$ of $Z_t$ given in \eqref{scale}, and of the speed
density $m:(0,1)\to\R$ of $Z_t$ given in \eqref{speed}. By
\cite[Eq.~(15.3.11)]{KarTay81} we get
\[
w(z_0) = 2\bigg(
u(z_0)
\int_{z_0}^b (S(b)-S(t))m(t) g(t) dt
+
u(1-z_0)
\int_a^{z_0} (S(t)-S(a))m(t) g(t) dt
\bigg),
\]
or, in terms of the \emph{Green function} $G :[a,b] \times [a,b]$ of
the process $Z_t$ on the interval $[a, b]$, 
\[
w(z_0) =\int_a^b G(x_0,t) g(t) dt,
\]
where
\[
G(x,s) = 
\begin{cases}
2\frac{
\int_{a}^{x} t^{-2 a_0 } (1-t)^{-2 a_1 } dt 
\int_{s}^{b} t^{-2 a_0 } (1-t)^{-2 a_1 } dt 
}{
\int_{a}^{b} t^{-2 a_0 } (1-t)^{-2 a_1 } dt 
} s^{2a_0-1}(1-s)^{2a_1-1}
& \text{if }a\leq x \leq s \leq b;
\\
2\frac{
\int_{x}^{b} t^{-2 a_0 } (1-t)^{-2 a_1 } dt 
\int_{a}^{s} t^{-2 a_0 } (1-t)^{-2 a_1 } dt 
}{
\int_{a}^{b} t^{-2 a_0 } (1-t)^{-2 a_1 } dt 
} s^{2a_0-1}(1-s)^{2a_1-1}
& \text{if }a\leq s \leq x \leq b.
\end{cases} 
\]
A complete characterization of $G$ in terms of a second order
differential equation may be found in \cite[p.~199]{KarTay81}.  An
explicit formula for $G$ is possible only when $a_0 = 0$ (extinction
of $J_1$, not admitted in our model), that may be found in
\cite[p.~208]{KarTay81}.

\subsubsection*{Problem 3.} 
In \cite{Grubel}, it is stated that for one-dimensional diffusion with
limiting invariant distribution with density $\pi(x)$, one has
\[
\pi(x) = \lim_{\varepsilon\to 0} \frac{E(\ind{\tau (x-\varepsilon,x+\varepsilon)})}
   {E(Z_{(x,\varepsilon)})},
\]
where
\[
\begin{aligned}
  \tau (x-\varepsilon,x+\varepsilon) & = \inf\{t>0 \colon Z_t
  \not\in (x-\varepsilon,x+\varepsilon)\},
\\
Z_{(x,\varepsilon)} & = \inf\{t>\tau(x-\varepsilon,x+\varepsilon)
\colon Z_t = x \},
\end{aligned}
\]
denote the first exit time from $(x-\varepsilon,x+\varepsilon)$ and
the time of first return to $x$ after leaving
$(x-\varepsilon,x+\varepsilon)$, respectively.  The proof in 
\cite{Grubel} can be easily modified to our context, so that for
$Z_t$ above relation reads
\begin{multline*}
  \lim_{\varepsilon\to 0} \frac{E(\ind{\tau (z_0-\varepsilon,z_0+\varepsilon)})}
      {E(Z_{(z_0,\varepsilon)})}
\\
= \frac{\Gamma (2 \frac{b}{\alpha} ) }{ \Gamma\Big( 2 \frac{b}{\alpha}
  \sum_{l\in J}{p_{l}} \Big) 
\Gamma\Big( 
2 \frac{b}{\alpha} (1-\sum_{l\in J}{p_{l}}) 
\Big) } 
z_0^{ 2 \frac{b}{\alpha} \sum_{l\in J}{p_{l}} -1 } (1-z_0)^{ 2 \frac{b}{\alpha}
  (1-\sum_{l\in J}{p_{l}}) -1}.
\end{multline*}


\section{Accessible and inaccessible boundaries: recessive sets and dominant components}
\label{dominant}

Looking at \eqref{eq:limit-distr}, we give the following definition:

\begin{Definition}
A subset $J \subsetneq \{1,\ldots,k\}$, $J\neq\varnothing$, is said
\emph{recessive} if $ \sum_{l \in J} p_{l} < \tfrac{\alpha}{2b} $.
\end{Definition}

Obviously, every subset of a recessive set is recessive. Moreover,
when $\tfrac{\alpha}{b} > 2(1-\min_i p_{i})$, every set $J \subsetneq
\{1,\ldots,k\}$ is recessive. Finally, the following result holds
true:

\begin{Proposition}\label{prop:recessive} We have:
\begin{enumerate}
\item\label{rec:p1} $J$ is recessive if and only if $P(\exists t
  \colon \cap_{i\in J} \{X_{t,i} = 0\} ) = 1$;
\item\label{rec:p2} $J$ is not recessive if and only if $P(\exists t
  \colon \cap_{i\in J} \{X_{t,i} = 0\} ) = 0$;
\item\label{rec:p3} the set $\{1,\ldots,i-1,i+1,\ldots,k\}$ is
  recessive if and only if $P(\exists t \colon \{X_{t,i} = 1\} ) = 1$;
\item\label{rec:p4} the set $\{1,\ldots,i-1,i+1,\ldots,k\}$ is not
  recessive if and only if $P(\exists t \colon \{X_{t,i} = 1\} ) = 0$.
\end{enumerate}
In case \ref{rec:p3}, the component $i$ is called \emph{dominant}.
\end{Proposition}

\begin{proof}
For $J \subsetneq \{1,\ldots,k\}$, $J\neq\varnothing$, let $\bs{J} =
\{J,J^c\}$.  Proposition~\ref{prop:2dimURN} implies that $Z_t =
\sum_{l\in J} X_{t,l}$ satisfies the following equation
\[
\begin{aligned}
d{Z}_{t} 
& =
-b\frac{{Z}_{t} -
  \sum_{l\in J}{p_{l}}}{\alpha} dt + 
  \sqrt{\max(0,{Z}_{t}(1-{Z}_{t}))} d {W}_t
  \\
  & = 
  \bigg(
-\frac{b}{\alpha}\Big(1-\sum_{l\in J}{p_{l}} \Big) Z_t 
 + \frac{b}{\alpha} \sum_{l\in J}{p_{l}}  (1-Z_t)
\bigg) dt + 
  \sqrt{\max(0,{Z}_{t}(1-{Z}_{t}))} d {W}_t
  ,
  \end{aligned}
\]
that is \eqref{eq:procBarrier} with $a_0 = \frac{b}{\alpha} \sum_{l\in
  J}{p_{l}} $ and $ a_1 = \frac{b}{\alpha} - a_0 $.  Since $\cap_{i\in
  J} \{X_{t,i} = 0\} = \{Z_t = 0\}$, the results \ref{rec:p1} and
\ref{rec:p2} in Proposition~\ref{prop:recessive} are a consequence of
the classification of the boundary point $z=0$ given in
Appendix~\ref{sec:bound}.

With the same spirit, Corollary~\ref{cor:margin} states that $Z_t = 1 - X_{t,i}$
satisfies the SDE
\[
\begin{aligned}
d{Z}_{t} 
& =
-b\frac{{Z}_{t} -
  \sum_{l\neq i} {p_{l}}}{\alpha} dt + 
  \sqrt{\max(0,(1-{Z}_{t}) {Z}_{t})} d {W}_t
  \\
  & = 
  \Big(
-\frac{b}{\alpha} p_i Z_t 
 + \frac{b}{\alpha} (1-p_i)  (1-Z_t)
\Big) dt + 
  \sqrt{\max(0,{Z}_{t}(1-{Z}_{t}))} d {W}_t
 ,
  \end{aligned}
\]
that is \eqref{eq:procBarrier} with $a_0 = \frac{b}{\alpha} (1-p_i)$ and $ a_1 = \frac{b}{\alpha}  - a_0 $.
The results \ref{rec:p3} and \ref{rec:p4} in Proposition~\ref{prop:recessive} are then a consequence
of the classification of the boundary point $z=0$ given in Appendix~\ref{sec:bound}.
\end{proof}

\appendix

\section{Multidimensional orthogonal polynomials on the simplex}
\label{sec:polynomials}
In this section, we recall some results on orthogonal polynomials on
the simplex, as given in \cite[Section~5.3]{Dunk2}.  Our notation
differs from that of \cite{Dunk2} since we use $\gamma_i$ instead of
$\kappa_i - 1/2$ and $k-1$ instead of $d$.  Accordingly, let $T^{k-1}$
be the $k-1$-dimensional simplex defined by
$$
T^{k-1} := \{\bs{y} \in \R^{k-1}: y_1 \ge 0, \ldots, y_{k-1} \ge 0,
1-y_1-y_2-\cdots-y_{k-1} \ge 0 \}. 
$$ Fixed $\bs{\gamma} = (\gamma_1,\ldots,\gamma_k)$ with $\gamma_i>-1$
for any $i\in\{1,\ldots,k\}$, the classical polynomials on $T^{k-1}$
are orthogonal with respect to the weight $\mathcal{L}^1(T^{k-1})$
function
\[
f_{\bs{\gamma}}(\bs{y}) =
(1-y_1-\cdots-y_{k-1})^{\gamma_k} \prod_{i=1}^{k-1} y_i^{\gamma_i} ,
\]
where the normalization constant $w_{\bs{\gamma}}$ of $f_{\bs{\gamma}}$ is
given by the Dirichlet integral
\[
\frac{1}{w_{\bs{\gamma}}} = \frac{ \Gamma(\gamma_1 +1) \cdots
  \Gamma(\gamma_k +1) }{ \Gamma \Big(\sum_1^k \gamma_l + k \Big)  }.
\]
Then, we may define $\pi_{\boldsymbol{\gamma}}(\bs{{y}}) =
w_{\bs{\gamma}} f_{\bs{\gamma}}(\bs{y}) $, which is a density on $T^{k-1}$.
The Hilbert space that we consider here is hence defined on $T^{k-1}$
by the inner product
\[
\langle f,g \rangle_{\boldsymbol{\gamma}} = \int_{T^{k-1}} f(\bs{y}) g(\bs{y}) 
\pi_{\boldsymbol{\gamma}}(\bs{{y}}) 
d\bs{y},
\]
that gives the orthogonality stated above.  As proven in
\cite[Section~5.3]{Dunk2}, the space
$\mathcal{V}_{n,\boldsymbol{\gamma}}$ of orthogonal polynomials of
degree $n$ is a eigenspace of eigenfunctions $f$ of the second-order
differential operator
\begin{equation*}
\begin{aligned}
L_{\bs{\gamma} } f & =
\sum_{i=1}^{k-1} \Big( \gamma_i +1 - \Big[ \big(\sum_{i=1}^k \gamma_i \big)
  + k \Big] y_i \Big) 
\frac{\partial f }{\partial y_i} (\bs{y} )
\\
& \qquad\qquad
+ \sum_{i=1}^{k-1} y_i (1 - y_i) \frac{\partial^2 f }{\partial y_i^2} (\bs{y}) 
- 2\sum_{1 \leq i < j \leq k-1} y_iy_j \frac{\partial^2 f }
{\partial y_i\partial y_j} (\bs{y}) 
\end{aligned}
\end{equation*}
with eigenvalue $\lambda_n = n(n+ k + \sum_{i=1}^k \gamma_i ) $ (see
\cite[Eq.~(5.3.4)]{Dunk2}), that is,
\begin{equation}\label{eq:diff-eqn2}
  L_{\boldsymbol{\gamma}} f (\bs{y})=  -\lambda_n f (\bs{y}) ,
  \qquad \text{for any } f \in \mathcal{V}_{n,\boldsymbol{\gamma}}.
\end{equation}
In \cite[Section~5.3]{Dunk2}), three orthogonal bases of
$\mathcal{V}_{n,\boldsymbol{\gamma}}$ are presented.  Each one of
these three bases is made by functions identified by the $
\binom{n+k-2}{n}$ possible choice of $\mathbf{n} =
(n_1,\ldots,n_{k-1})$ with $n_i \geq 0$ and $|\mathbf{n}| = \sum_i n_i
= n$ as follows.
\begin{description}
\item[Jacobi:] the standard extension of Jacobi polynomials given in
  \cite[Proposition~5.3.1]{Dunk2} as
\[
 P^{\boldsymbol{\gamma}}_{{\mathbf{n}}} ({\bs{y}}) = h^{\boldsymbol{\gamma}}_{P,{\mathbf{n}}}
 \prod_{i=1}^{k-1} \Big(\sum_{l=i}^{k-1} x_l \Big)^{n_i} 
 \mathrm{p}_{n_i}^{(\sum_{l=i+1}^{k-1}n_l + \sum_{l=i+1}^{k} (\gamma_l+1) - 1,\gamma_i)}
 \Big(\tfrac{2 x_i}{\sum_{l=i}^{k-1} x_l }-1\Big),
\]
where $\mathrm{p}_{n}^{(a_1,a_2)} (t) $ is the standard Jacobi
polynomial defined on $-1<t<1$ and
$h^{\boldsymbol{\gamma}}_{P,{\mathbf{n}}}$ is the normalizing factor
given in \cite[Proposition~5.3.1]{Dunk2} in terms of products of
Pochhammer symbols, so that $\langle
P^{\boldsymbol{\gamma}}_{{\mathbf{n_1}}} ,
P^{\boldsymbol{\gamma}}_{{\mathbf{n_2}}} \rangle_{\boldsymbol{\gamma}}
= \ind{{\mathbf{n_1}}\equiv {\mathbf{n_2}}} $;
\item[Monic orthogonal basis:] in \cite[Proposition~5.3.2]{Dunk2} it
  is proven that the following family is a orthogonal base of
  $\mathcal{V}_{n,\boldsymbol{\gamma}}$:
\[
V^{\boldsymbol{\gamma}}_{{\mathbf{n}}} ({\bs{y}}) =
 \sum_{\mathbf{m} \preceq \mathbf{n}} 
 \prod_{i=1}^{k-1} 
 (-1)^{n_i-m_i} y_i^{m_i}
 \bigg(
 \binom{n_i}{m_i} \frac{(\gamma_i+1)_{n_i}}{(\gamma_i+1)_{m_i}}
 \frac{(\sum_1^k \gamma_l + k - 1)_{n_i}}{(\sum_1^k \gamma_l + k - 1)_{m_i}}
 \bigg),
 \]
 where $\mathbf{m} \preceq \mathbf{n}$ means $m_i \leq n_i $ for any
 $i = 1,\ldots,k-1$. Note that in this case the normalizing factor
 $h^{\boldsymbol{\gamma}}_{V,{\mathbf{n}}}$ is not given explicitly,
 and then $\langle V^{\boldsymbol{\gamma}}_{{\mathbf{n_1}}} ,
 V^{\boldsymbol{\gamma}}_{{\mathbf{n_2}}}
 \rangle_{\boldsymbol{\gamma}} = \ind{{\mathbf{n_1}}\equiv
   {\mathbf{n_2}}} (h^{\boldsymbol{\gamma}}_{V,{\mathbf{n_1}}})^2 $.
\item[Rodrigue formula:] in \cite[Proposition~5.3.3]{Dunk2} it is
  proven that the following family is a orthogonal base of
  $\mathcal{V}_{n,\boldsymbol{\gamma}}$, given in terms of the
  Rodrigue formula:
\begin{equation*} 
 U^{\boldsymbol{\gamma}}_{{\mathbf{n}}} ({\bs{y}}) = 
 \frac{
 \frac{\partial^{n}} {\partial y_1^{n_1} \cdots 
\partial y_{k-1}^{n_{k-1}}}
     \left[ 
     (1- y_1-y_2-\cdots-y_{k-1} )^{\gamma_{k} + n} \prod_{i=1}^{k-1}y_i^{\gamma_i+n_i} 
     \right]
 }{
 (1- y_1-y_2-\cdots-y_{k-1} )^{\gamma_{k}} \prod_{i=1}^{k-1}y_i^{\gamma_i} 
}.
\end{equation*}
Again, the normalizing factor
$h^{\boldsymbol{\gamma}}_{U,{\mathbf{n}}}$ is not given explicitly,
and then $\langle U^{\boldsymbol{\gamma}}_{{\mathbf{n_1}}} ,
U^{\boldsymbol{\gamma}}_{{\mathbf{n_2}}} \rangle_{\boldsymbol{\gamma}}
= \ind{{\mathbf{n_1}}\equiv {\mathbf{n_2}}}
(h^{\boldsymbol{\gamma}}_{U,{\mathbf{n_1}}})^2 $. Moreover, the two
families $(U^{\boldsymbol{\gamma}}_{{\mathbf{n}}} )_{{\mathbf{n}}}$
and $(V^{\boldsymbol{\gamma}}_{{\mathbf{n}}} )_{{\mathbf{n}}}$ are
biorthogonal, in the sense that $\langle
U^{\boldsymbol{\gamma}}_{{\mathbf{n_1}}} ,
V^{\boldsymbol{\gamma}}_{{\mathbf{n_2}}} \rangle_{\boldsymbol{\gamma}}
= 0$ whenever ${\mathbf{n_1}} \neq {\mathbf{n_2}}$.
\end{description}


\section{Wright-Fisher boundary types}\label{sec:bound}
In this section we recall a classification of the boundaries of the 
one-dimensional Wright-Fisher process with mutation given in
\cite[p.~239, Example~8]{KarTay81} (see also \cite{Huil07}).

Fixed $a_0,a_1 \geq 0$, let ${Z}_{t}$ be the process with values in
$[0,1]$ that satisfies the SDE
\begin{equation}\label{eq:procBarrier}
d{Z}_{t} = (-a_1 Z_t + a_0 (1-Z_t)) dt
 + 
  \sqrt{\max(0,{Z}_{t}(1-{Z}_{t}))} d {W}_t.
\end{equation}
Then by \cite[Eq.~(15.6.18) and Eq.~(15.6.19)]{KarTay81}, we have that
\[
\text{the boundary point }z\in \{0,1\}
\begin{cases}
\text{is an exit boundary} & \text{if }a_z = 0;
\\
\text{is a regular boundary} & \text{if } 0 < a_z < \tfrac{1}{2};
\\
\text{is an entrance boundary} & \text{if }a_z \geq \tfrac{1}{2}.
\end{cases} 
\]
When $a_0,a_1>0$, the SDE may be written as
\[
d{Z}_{t} 
= -(a_0+a_1)\Big(Z_t - \frac{a_0}{a_0+a_1} \Big) dt
 + 
  \sqrt{\max(0,{Z}_{t}(1-{Z}_{t}))} d {W}_t  
\]
and the process that starts at $Z_0\in(0,1)$ will never leave the
strip $[0,1]$. In particular, the classification above states whether 
the process $Z_t$ will reach the boundary infinitely many times (and
the boundary point is a reflection barrier) or will never reach it. In
particular $z$ is a regular boundary if and only if $Z_t$ reaches the
reflecting barrier $z$ infinitely many times; while $z$ is an entrance
boundary if and only if the process $Z_t$ will never touch $z$.

Moreover, we may compute the \emph{scale function} $S:(0,1)\to\R$,
defined as the integral of its derivative
\[
S'(z) = S'(z_0)
\exp \Big( - \int_{z_0}^z \frac{2 (-a_1 t + a_0 (1-t))}{t(1-t)} dt \Big)
\]
and the \emph{speed density} 
$m:(0,1)\to\R$, defined as $ m(z) = 1/(S'(z) z (1-z))$.
A direct computation as in \cite{Huil07} yields:
\begin{align}
  S(z) & = S(z_1) + S'(z_1) \int_{z_0}^z t^{-2 a_0} (1-t)^{-2a_1} dt
  \qquad\mbox{and} \label{scale}
  \\
m(z) & = \tfrac{1}{S'(z_1)} z^{2a_0-1} (1-z)^{2a_1-1}  \label{speed}.
\end{align}


\begin{thebibliography}{10}

\bibitem{ale-cri-GRP}
G.~Aletti and I.~Crimaldi.
\newblock Generalized rescaled {P}\'olya urn and its statistical applications.
\newblock {\em arXiv2010.06373}, 2021.

\bibitem{ale-cri-RP}
G.~Aletti and I.~Crimaldi.
\newblock The rescaled {P}\'olya urn: local reinforcement and chi-squared
  goodness of fit test.
\newblock {\em Advances in Applied Probability}, 54:{forthcoming}, 2022.

\bibitem{ale-cri-sar-sentiment}
G.~Aletti, I.~Crimaldi, and F.~Saracco.
\newblock A model for the twitter sentiment curve.
\newblock {\em PLOS ONE}, 16(4):1--28, 04 2021.

\bibitem{bol-et-al}
J.~P. Bollback, T.~L. York, and R.~Nielsen.
\newblock {Estimation of $2N_es$ from temporal allele frequency data}.
\newblock {\em Genetics}, 179:497--502, 2008.

\bibitem{cha-gen}
M.~Chaleyat-Maurel and V.~Genon-Catalot.
\newblock Filtering the {W}right{\textendash}{F}isher diffusion.
\newblock {\em ESAIM: Probability and Statistics}, 13:197--217, 6 2009.

\bibitem{cost-et-al}
C.~Costantini, P.~De~Blasi, S.~Ethier, M.~Ruggiero, and D.~Span\`o.
\newblock {W}right-{F}isher construction of the two-parameter
  {P}oisson-{D}irichlet diffusion.
\newblock {\em The Annals of Applied Probability}, 27:1923–1950, 2017.

\bibitem{dan-et-al}
C.~Dangerfield, D.~Kay, and K.~Burrage.
\newblock Stochastic models and simulation of ion channel dynamics.
\newblock {\em Procedia Computer Science}, 1(1):1587--1596, 2010.

\bibitem{dan-et-al-2}
C.~E. Dangerfield, D.~Kay, S.~MacNamara, and K.~Burrage.
\newblock A boundary preserving numerical algorithm for the
  {W}right{\textendash}{F}isher model with mutation.
\newblock {\em {BIT Numerical Mathematics}}, 5:283--304, 2012.

\bibitem{del-shi}
F.~Delbaen and H.~Shirakawa.
\newblock An interest rate model with upper and lower bounds.
\newblock {\em Asia-Pac. Financ. Mark.}, 9:191--209, 2002.

\bibitem{Dunk2}
C.~F. Dunkl and Y.~Xu.
\newblock {\em Orthogonal polynomials of several variables}, volume 155 of {\em
  Encyclopedia of Mathematics and its Applications}.
\newblock Cambridge University Press, Cambridge, second edition, 2014.

\bibitem{EggPol23}
F.~Eggenberger and G.~P\'olya.
\newblock {\"{U}}ber die statistik verketteter vorg\"ange.
\newblock {\em ZAMM - Journal of Applied Mathematics and Mechanics /
  Zeitschrift f\"ur Angewandte Mathematik und Mechanik}, 3(4):279--289, 1923.

\bibitem{fav-et-al}
S.~Favaro, M.~Ruggiero, and S.~G. Walker.
\newblock On a {G}ibbs sampler based random process in {B}ayesian
  nonparametrics.
\newblock {\em Electronic Journal of Statistics}, 3:1556--1566, 2009.

\bibitem{gou-jas}
C.~Gourieroux and J.~Jasiak.
\newblock Multivariate jacobi process with application to smooth transitions.
\newblock {\em Journal of Econometrics}, 131:475--505, 2006.

\bibitem{gri-spa}
R.~C. Griffiths and D.~Span\`o.
\newblock Diffusion processes and coalescent trees.
\newblock {\em In Probability and Mathematical Genetics, Papers in Honour of
  Sir John Kingman, (N. H. Bingham and C. M. Goldie, eds.). LMS Lecture Note
  Series. Cambridge University Press}, 378(15):358--375, 2010.

\bibitem{Grubel}
R.~Gr\"{u}bel.
\newblock On mean recurrence times of stationary one-dimensional diffusion
  processes.
\newblock {\em Stochastic Process. Appl.}, 18(1):165--169, 1984.

\bibitem{gut-et-al}
R.~N. Gutenkunst, R.~D. Hernandez, S.~H. Williamson, and C.~D. Bustamante.
\newblock Inferring the joint demographic history of multiple populations from
  multidimensional snp frequency data.
\newblock {\em PLOS Genetics}, 5(10):1--11, 10 2009.

\bibitem{Huil07}
T.~Huillet.
\newblock On {W}right{\textendash}{F}isher diffusion and its relatives.
\newblock {\em Journal of Statistical Mechanics: Theory and Experiment},
  2007(11):P11006--P11006, nov 2007.

\bibitem{jen-spa}
P.~A. Jenkins and D.~Spanò.
\newblock Exact simulation of the {W}right-{F}isher diffusion.
\newblock {\em The Annals of Applied Probability}, 27(3):1478 -- 1509, 2017.

\bibitem{KarTay81}
S.~Karlin and H.~M. Taylor.
\newblock {\em A second course in stochastic processes}.
\newblock Academic Press, Inc. [Harcourt Brace Jovanovich, Publishers], New
  York-London, 1981.

\bibitem{Kus}
H.~J. Kushner.
\newblock {\em Approximation and weak convergence methods for random processes,
  with applications to stochastic systems theory}, volume~6 of {\em MIT Press
  Series in Signal Processing, Optimization, and Control}.
\newblock MIT Press, Cambridge, MA, 1984.

\bibitem{KusYin}
H.~J. Kushner and G.~G. Yin.
\newblock {\em Stochastic approximation and recursive algorithms and
  applications}, volume~35 of {\em Applications of Mathematics (New York)}.
\newblock Springer-Verlag, New York, second edition, 2003.
\newblock Stochastic Modelling and Applied Probability.

\bibitem{mah}
H.~M. Mahmoud.
\newblock {\em {P}\'olya urn models}.
\newblock Texts in Statistical Science Series. CRC Press, Boca Raton, FL, 2009.

\bibitem{mala-et-al}
A.~S. Malaspinas, O.~Malaspinas, S.~N. Evans, and M.~Slatkin.
\newblock Estimating allele age and selection coefficient from time-serial
  data.
\newblock {\em Genetics}, 192:599--607, 2012.

\bibitem{mena-rug}
R.~Mena and M.~Ruggiero.
\newblock Dynamic density estimation with diffusive dirichlet mixtures.
\newblock {\em Bernoulli}, 22:901--926, 2016.

\bibitem{pap-rug}
O.~Papaspiliopoulos and M.~Ruggiero.
\newblock {Optimal filtering and the dual process}.
\newblock {\em Bernoulli}, 20(4):1999 -- 2019, 2014.

\bibitem{sch-et-al}
J.~Schraiber, R.~C. Griffiths, and S.~N. Evans.
\newblock Analysis and rejection sampling of {W}right-{F}isher diffusion
  bridges.
\newblock {\em Theoretical Population Biology}, 89:64--74, 2013.

\bibitem{Tanabe}
K.~Tanabe and M.~Sagae.
\newblock An exact cholesky decomposition and the generalized inverse of the
  variance-covariance matrix of the multinomial distribution, with
  applications.
\newblock {\em Journal of the Royal Statistical Society. Series B
  (Methodological)}, 54(1):211--219, 1992.

\bibitem{wal-et-al}
S.~G. Walker, S.~J. Hatjispyros, and T.~Nicoleris.
\newblock A {F}leming-{V}iot process and {B}ayesian nonparametrics.
\newblock {\em Annals of Applied Probability}, 17:67--80, 2007.

\bibitem{wil-et-al}
S.~H. Williamson, R.~Hernandez, A.~Fledel-Alon, L.~Zhu, and C.~D. Bustamante.
\newblock Simultaneous inference of selection and population growth from
  patterns of variation in the human genome.
\newblock {\em Proceedings of the National Academy of Sciences of the United
  States of America}, 102:7882--7887, 2005.

\bibitem{wright84}
S.~Wright.
\newblock {\em Evolution and the Genetics of Populations, Volume 2: Theory of
  Gene Frequencies}.
\newblock Evolution and the Genetics of Populations. University of Chicago
  Press, 1984.

\bibitem{zhao-et-al}
L.~Zhao, M.~Lascoux, A.~D.~J. Overall, and D.~Waxman.
\newblock The characteristic trajectory of a fixing allele: a consequence of
  fictitious selection that arises from conditioning.
\newblock {\em Genetics}, 195:993--1006, 2013.

\end{thebibliography}

\section*{Acknowledgements}
Both authors sincerely thank Fabio Saracco for having collected and shared with them the two Twitter data sets. Giacomo Aletti is a member of the Italian Group ``Gruppo
Nazionale per il Calcolo Scientifico'' of the Italian Institute
``Istituto Nazionale di Alta Matematica'' and Irene Crimaldi is a
member of the Italian Group ``Gruppo Nazionale per l'Analisi
Matematica, la Probabilit\`a e le loro Applicazioni'' of the Italian
Institute ``Istituto Nazionale di Alta Matematica''.  \\[5pt]
\noindent{\bf Funding Sources}\\
\noindent Irene Crimaldi is partially supported by the Italian
``Programma di Attivit\`a Integrata'' (PAI), project ``TOol for
Fighting FakEs'' (TOFFE) funded by IMT School for Advanced Studies
Lucca.

\noindent{\bf Author contributions statement}\\
Both authors equally contributed to this work.

\end{document}